# Solutions of Linear Fractional non-Homogeneous Differential Equations with Jumarie Fractional Derivative and Evaluation of Particular Integrals

**Uttam Ghosh[1], Susmita Sarkar[2] and Shantanu Das[3]**

[1] Department of Mathematics, Nabadwip Vidyasagar College, Nabadwip, Nadia, West Bengal, India;
email : uttam_math@yahoo.co.in

[2] Department of Applied Mathematics, University of Calcutta, Kolkata, India
email : susmita62@yahoo.co.in

[3] Scientist H+, Reactor Control System Design Section (E & I Group) BARC Mumbai India
Senior Research Professor, Dept. of Physics, Jadavpur University Kolkata
Adjunct Professor. DIAT-Pune
UGC Visiting Fellow. Dept of Appl. Mathematics; Univ. of Calcutta
email : shantanu@barc.gov.in

## Abstract

In this paper we describe a method to solve the linear non-homogeneous fractional differential equations (FDE), composed with Jumarie type Fractional Derivative, and describe this method developed by us, to find out Particular Integrals, for several types of forcing functions. The solutions are obtained in terms of Mittag-Leffler functions, fractional sine and cosine functions. We have used our earlier developed method of finding solution to homogeneous FDE composed via Jumarie fractional derivative, and extended this to non-homogeneous FDE. We have demonstrated these developed methods with few examples of FDE, and also applied in fractional damped forced differential equation. This method proposed by us is useful as it is having conjugation with the classical methods of solving non-homogeneous linear differential equations, and also useful in understanding physical systems described by FDE.

## 1.0: Introduction

The fractional differential equations and its solutions arises in different branches of applied science, engineering, applied mathematics and biology [1-9]. The solutions of fractional difference equations are obtained by different methods which includes Exponential-Function Method [10], Homotopy Perturbation Method [11], Variation Iteration Method [12], Differential transform Method [13] and Fractional Sub-equation Method [14], Analytical Solutions in terms of Mittag-Leffler function [15]. In developing those Methods the usually used fractional derivative is Riemann-Liouvellie (R-L) [6], Caputo derivative [6], Jumarie's left handed modification of R-L fractional derivative [16-17]. In [15] we have developed an algorithm to solve the homogeneous fractional order differential equations in terms of Mittag-Leffler function and fractional sine and cosine functions. However, there are no standard methods to find solutions of non-homogeneous fractional differential equations. In this paper we describe a method to solve the fractional order non-homogeneous differential equations. Organizations of the paper are as follows; in section 2.0 we describe the different definitions of fractional derivatives and properties of Mittag-Leffler function. In section-3.0 we describe the solutions of $\alpha$ – order fractional differential equations. In section 4.0 the solutions of $2\alpha$ – order fractional differential equations is described, with several types of forcing functions. In section 5.0 this methods has been applied to solve both un-damped and damped fractional order forced oscillator equations. In this paper the fractional derivative operator $D^{\upsilon}$ will be of Jumarie type fractional derivative.

## 2.0: Definition of fractional derivatives

The useful definitions of the fractional derivatives are the Grunwald-Letinikov (G-L) definition and Riemann-Liouville(R-L) definition [6] and Modified R-L-definitions [16-17].

- **Grunwald-Letinikov definition**

Let $f(t)$ be any function then the $\alpha$-th order derivative $\alpha \in \mathbb{R}$ of $f(t)$ is defined by

$$_aD_t^{\alpha}[f(t)] = \lim_{\substack{h \to 0 \\ nh \to t-a}} h^{-\alpha} \sum_{r=0}^{n} \binom{\alpha}{r} f(t-rh) = \lim_{\substack{h \to 0 \\ nh \to t-a}} h^{-\alpha} \sum_{r=0}^{n} \frac{\alpha!}{r!(\alpha-r)!} f(t-rh)$$

$$= \frac{1}{(-\alpha-1)!} \int_a^t (t-\tau)^{-\alpha-1} f(\tau) d\tau = \frac{1}{\Gamma(-\alpha)} \int_a^t \frac{f(\tau)}{(t-\tau)^{\alpha+1}} d\tau$$

Where $\alpha$ is any arbitrary number real or complex; and the generalized binomial coefficients are described as follows [1], [16-17]

$$\binom{\alpha}{r} = \frac{\alpha!}{r!(\alpha-r)!} = \frac{\Gamma(\alpha+1)}{\Gamma(r+1)\Gamma(\alpha-r+1)}$$

The above formula becomes fractional order integration if we replace $\alpha$ by $-\alpha$ which is Riemann fractional integration formula

$$\int_a^t f(\tau) d\tau^{\alpha} = {_aI_t^{\alpha}}[f(t)] = {_aD_t^{-\alpha}}[f(t)] = \frac{1}{\Gamma(\alpha)} \int_a^t (t-\tau)^{\alpha-1} f(\tau) d\tau \qquad \alpha > 0$$

In above we have noted several notations used for fractional integration.

- **Riemann-Liouville fractional derivative definition**

Let the function $f(t)$ is one time integrable then the integro-differential expression as following defines Riemann-Liouvelli fractional derivative [1], [6]

$$_aD_t^{\alpha}[f(t)] = \frac{1}{\Gamma(n-\alpha)} \left(\frac{d}{dt}\right)^n \int_a^t (t-\tau)^{(n-\alpha)-1} f(\tau) d\tau \qquad \alpha > 0$$

Here the $n$ is a positive integer number just greater than real number $\alpha$ The above expression is known as the Riemann-Liouville definition of fractional derivative [6] with $(n-1) < \alpha < n$

In the above definition fractional derivative of a constant is non-zero.

- **Modified Riemann-Liouville definition**

To overcome the shortcoming fractional derivative of a constant, as non-zero, another modification of the definition of left R-L type fractional derivative of the function $f(x)$, in the interval $[a,b]$ was proposed by Jumarie [16] in the form described below

$$_a^J D_x^\alpha [f(x)] = \begin{cases} \dfrac{1}{\Gamma(-\alpha)} \int_a^x (x-\tau)^{-\alpha-1} f(\tau) d\tau, & \alpha < 0. \\ \dfrac{1}{\Gamma(1-\alpha)} \dfrac{d}{dx} \int_a^x (x-\tau)^{-\alpha} [f(\tau) - f(a)] d\tau, & 0 < \alpha < 1 \\ \left( f^{(\alpha-m)}(x) \right)^{(m)}, & m \leq \alpha < m+1. \end{cases}$$

Here we state that $f(x) = 0$ for $x < a$ and $x > b$. However in this paper we will be using this left-Jumarie fractional derivative that is $_0^J D_t^\alpha [f(t)]$, for $0 < \alpha < 1$ and with condition $f(t) = 0$ for all $t < 0$. We will simplify the symbol and drop $a = 0$ and differentiationg variable $t$ and simply write $^J D^\alpha [f(t)]$. Using the above definition Jumarie [16] proved

$$(u(x)v(x))^{(\alpha)} = u^{(\alpha)}(x)v(x) + u(x)v^{(\alpha)}(x)$$

We have recently modified the right R-L definition of fractional derivative of the function $f(x)$, in the interval $[a,b]$ in the following form [17],

$$_x^J D_b^\alpha [f(x)] = \begin{cases} -\dfrac{1}{\Gamma(-\alpha)} \int_x^b (\tau-x)^{-\alpha-1} f(\tau) d\tau, & \alpha < 0. \\ -\dfrac{1}{\Gamma(1-\alpha)} \dfrac{d}{dx} \int_x^b (\tau-x)^{-\alpha} [f(b) - f(\tau)] d\tau, & 0 < \alpha < 1 \\ \left( f^{(\alpha-m)}(x) \right)^{(m)}, & m \leq \alpha < m+1. \end{cases}$$

Using both the modified definition we investigate the characteristics of the non-differentiable points of some continuous functions in [17]. The above defined all the derivatives are non-local type, and obtained solution to homogeneous FDE, with Jumarie derivative [15]. Subsequently we will be using $^J D^v$ as fractional derivative operator of Jumarie type, with start point $a = 0$, and stating the function $f(x) = 0$ for all $x < 0$ in following sections.

## 2.1 The Mittag-Leffler Function

The Mittag-Leffler function was introduced by Gösta Mittag-Leffler [18] in 1903. The one-parameter Mittag-Leffler function is denoted by $E_\alpha(t^\alpha)$ and defined by following series

$$E_\alpha(z) \stackrel{\text{def}}{=} \sum_{k=0}^{\infty} \frac{z^k}{\Gamma(1+\alpha k)}, \qquad z \in \mathbb{C}, \qquad \text{Re}(\alpha) > 0$$

Again from the Jumarie definition of fractional derivative we have ${}_a^J D_x^\alpha [C] = 0$ we apply this property to get $\alpha$ order Jumarie Derivative of the Mittag-Leffler function $E_\alpha(ax^\alpha)$ as follows

$$\begin{aligned}
{}_0^J D_x^\alpha \left[ E_\alpha(ax^\alpha) \right] &= {}^J D^\alpha \left[ E_\alpha(ax^\alpha) \right] \\
&= {}^J D^\alpha \left[ 1 + \frac{ax^\alpha}{\Gamma(1+\alpha)} + \frac{a^2 x^{2\alpha}}{\Gamma(1+2\alpha)} + \frac{a^3 x^{3\alpha}}{\Gamma(1+3\alpha)} + \ldots \infty \right] \\
&= a\left( 1 + \frac{ax^\alpha}{\Gamma(1+\alpha)} + \frac{a^2 x^{2\alpha}}{\Gamma(1+2\alpha)} + \frac{a^3 x^{3\alpha}}{\Gamma(1+3\alpha)} + \ldots \infty \right) \\
&= aE_\alpha(ax^\alpha)
\end{aligned}$$

Therefore the fractional differential equation ${}_0^J D_x^\alpha [y] = ay$ has solution in the form $y = AE_\alpha(ax^\alpha)$, where $A$ is an arbitrary constant.

## 2.2 Non-Homogeneous Fractional Differential Equations and Some Basic Solutions

The general format of the fractional linear differential equation is

$$f\left( {}^J D^\alpha \right) y = g(t^\alpha) \tag{2.1}$$

Where $f\left( {}^J D^\alpha \right)$ is a linear differential operator $0 < \alpha < 1$. The above differential equation is said to be linear non-homogeneous fractional differential equation when $g(t^\alpha) \neq 0$, otherwise it is homogeneous. Solution of the linear fractional differential equations (composed via Jumarie Derivative) can be easily obtained in terms of Mittag-Leffler function and fractional sine and cosine functions [15].

The function $g(t^\alpha)$ is forcing function. We have written this as function of $t^\alpha$ purposely for ease. For example we will use in this paper $E_\alpha(ct^\alpha)$, $\sin_\alpha(ct^\alpha)$, $t^\alpha$ etc. are taken as forcing functions. There will be other functions in the derivations like $V(t^\alpha)$, $V(t^\alpha)$ all functions described with scaled variable that is $t^\alpha$. Nevertheless the forcing functions can be written as simple $h(t)$ though.

In that paper [15] we found the following (theorems) which we will be using in this paper

(i)  The fractional differential equation $\left({}_0^J D_t^\alpha - a\right)\left({}_0^J D_t^\alpha - b\right) y(t) = 0$ has solution of the form $y = AE_\alpha(at^\alpha) + BE_\alpha(bt^\alpha)$ where A and B are constants,

(ii) The fractional differential equation ${}_0^J D_t^{2\alpha}[y] - 2a\left({}_0^J D_t^\alpha [y]\right) + a^2 y = 0$ has solution of the form $y = (At^\alpha + B)E_\alpha(at^\alpha)$ where A and B are constants and

(iii) Solution of the fractional differential equation ${}_0^J D_t^{2\alpha}[y] - 2a\left({}_0^J D_t^\alpha [y]\right) + (a^2 + b^2) y = 0$ is of the form $y = E_\alpha(at^\alpha)[A\cos_\alpha(bt^\alpha) + B\sin_\alpha(bt^\alpha)]$ where A and B are constants.

From now we indicate Jumarie fractional derivative with start point of differentiation as 0 as ${}^J D^\alpha$ instead ${}_0^J D^\alpha$.

**Theorem 1:** If $y_1$ and $y_2$ are two solutions of the fractional differential equation $f\left({}^J D^\alpha\right) y = 0$ then $c_1 y_1 + c_2 y_2$ is also a solution, where $c_1$ and $c_2$ are arbitrary constants.

**Proof:** Since $f(D^\alpha) y = 0$ has solutions $y = y_1$ and $y = y_2$

$$f({}^J D^\alpha) y_1 = 0 \quad \text{and} \quad f({}^J D^\alpha) y_2 = 0$$
$$f({}^J D^\alpha)(c_1 y_1 + c_2 y_2) = c_1 f({}^J D^\alpha) y_1 + c_2 f({}^J D^\alpha) y_2 = 0$$

Hence $c_1 y_1 + c_2 y_2$ is also a solution of the given fractional differential equation.

Hence the theorem is proved.

Similarly, we can prove if $y_1, y_2, ..., y_n$ are solutions of the fractional differential equation $f({}^J D^\alpha) y = 0$ then $c_1 y_1 + c_2 y_2 + ... + c_n y_n$ is also a solution of it.

**Theorem 2:** If $f({}^J D^\alpha) = ({}^J D^\alpha - a_1)({}^J D^\alpha - a_2)({}^J D^\alpha - a_3)...({}^J D^\alpha - a_n)$, $0 < \alpha \le 1$. then solution of the homogeneous equation $f({}^J D^\alpha) y = 0$ is $y = \sum_{k=1}^{n} A_k E_\alpha(a_k t^\alpha)$ where $A_k$'s are arbitrary constants and all $a_k$ are distinct.

**Proof:** Since Jumarie type fractional derivative of Mittag-Leffler function $y = E_\alpha(at^\alpha)$ with $a$ as a constant is ${}^J D^\alpha[y] = aE_\alpha(at^\alpha) = ay$, $0 < \alpha \le 1$. Thus solution of the differential equation ${}^J D^\alpha y = ay$, $0 < \alpha \le 1$ is $y = AE_\alpha(at^\alpha)$ where A is a constant [15].

Let $y = AE_\alpha(mt^\alpha) \ne 0$ be a non-trivial trial solution of the differential equation $f({}^J D^\alpha) y = 0$ then ${}^J D^\alpha y = my$ or we write the following after subtracting $a_i y$ from both the sides as demonstrated below

$$^J D^\alpha y = {}^J D^\alpha \left[ A E_\alpha (m t^\alpha) \right] = A \left( m E_\alpha (m t^\alpha) \right) = m y$$
$$^J D^\alpha y - a_i y = m y - a_i y$$
$$({}^J D^\alpha - a_i) y = (m - a_i) y$$

We apply the above result sequentially as demonstrated below

$$f({}^J D^\alpha) y = f({}^J D^\alpha) \left[ A E_\alpha (m t^\alpha) \right]$$
$$= ({}^J D^\alpha - a_1)({}^J D^\alpha - a_2)\ldots\{({}^J D^\alpha - a_n)\left[ A E_\alpha (m t^\alpha) \right]\}$$
$$= ({}^J D^\alpha - a_1)({}^J D^\alpha - a_2)\ldots({}^J D^\alpha - a_{n-1})\{(m - a_n)\left[ A E_\alpha (m t^\alpha) \right]\}$$
$$= ({}^J D^\alpha - a_1)\{(m - a_2)(m - a_3)\ldots(m - a_n)\left[ A E_\alpha (m t^\alpha) \right]\}$$
$$= (m - a_1)(m - a_2)\ldots(m - a_n)\left[ A E_\alpha (m t^\alpha) \right] = \prod_{i=1}^{m}(m - a_i) y$$

Since $f({}^J D^\alpha) y = 0$ we get

$$\prod_{i=1}^{n}(m - a_i) y = 0 \qquad (2.2)$$

Implying that $m = a_1, a_2, \ldots a_n$

Hence the general solution is

$$y = A_1 E_\alpha (a_1 t^\alpha) + A_2 E_\alpha (a_2 t^\alpha) + \ldots + A_n E_\alpha (a_n t^\alpha) = \sum_{k=1}^{n} A_k E(a_k t^\alpha)$$

Hence the theorem is proved.

The above theorem implies principal of superposition holds for the linear fractional differential equations (composed via Jumarie fractional derivative) also.

Note: In the above theorem if two or more roots of the equation (2.2) are equal or roots are complex then the solution [15] form is given below.

For $a_1 = a_2$ and $a_3 \neq a_4 \neq \ldots \neq a_n$ then solution of the is

$$y = [A_1 t^\alpha + A_2] E_\alpha (a_1 t^\alpha) + A_3 E_\alpha (a_3 t^\alpha) + \ldots + A_n E_\alpha (a_n t^\alpha)$$

For $a_1 = a_2 = a_3$ and $a_4 \neq a_5 \neq \ldots \neq a_n$

then the solution is

$$y = [A_1 t^{2\alpha} + A_2 t^\alpha + A_3] E_\alpha (a_1 t^\alpha) + A_4 E_\alpha (a_4 t^\alpha) + \ldots + A_n E_\alpha (a_n t^\alpha)$$

where $A_k$'s are arbitrary constants.

For $a_1, a_2 = a \pm ib$ and other are $a_3 \neq a_4 \neq .... \neq a_n$ then the solution is

$$y = [A_1 \cos_\alpha(bt^\alpha) + A_2 \sin_\alpha(bt^\alpha)]E_\alpha(at^\alpha) + A_3 E_\alpha(a_3 t^\alpha) + ... + A_n E_\alpha(a_n t^\alpha)$$

Thus solutions of linear homogeneous fractional differential equation with Jumarie fractional derivative is express in terms of Mittag-Leffler functions and fractional type sine and cosine series.

Now the question arises what will be solution of linear non-homogeneous fractional differential equations. The solution corresponding to the homogeneous equation will be called as the complementary function, it contains the arbitrary constants and this solution will be denoted by $y_c$. The other part, that is a solution which is free from integral constant, and depending on the forcing function will be called as Particular Integral (PI) and will be denoted by $y_p$. Thus the general solution will be $y = y_c + y_p$. We will develop simple method to evaluate Particular Integral.

## 3.0 $\alpha$-order non-homogeneous fractional differential equations

Consider the linear $\alpha$- order non-homogeneous fractional differential equation with $0 < \alpha < 1$ for $y = 0$ for $t < 0$ of the following form,

$$(^J D^\alpha - a)y = g(t^\alpha) \qquad 0 < \alpha < 1 \qquad (3.1)$$

The solution of the corresponding homogeneous part is [15]

$$y_c = A_1 E_\alpha(mt^\alpha), \qquad A_1 \text{ is arbitary constant.}$$

Multiply both side of equation (3.1) by $E_\alpha(-at^\alpha)$ as demonstrated below

$$\left[E_\alpha(-at^\alpha)\right]\left[(^J D^\alpha - a)y\right] = g(t^\alpha)E_\alpha(-at^\alpha)$$
$$\left[E_\alpha(-at^\alpha)\right]\left[^J D^\alpha y - ay\right] = g(t^\alpha)E_\alpha(-at^\alpha)$$
$$\left[E_\alpha(-at^\alpha)\right]\left(^J D^\alpha [y]\right) - ayE_\alpha(-at^\alpha) = g(t^\alpha)E_\alpha(-at^\alpha)$$
$$\left[E_\alpha(-at^\alpha)\right]\left(^J D^\alpha [y]\right) + y\left(^J D^\alpha \left[E_\alpha(-at^\alpha)\right]\right) = g(t^\alpha)E_\alpha(-at^\alpha)$$
$$^J D^\alpha \left[yE_\alpha(-at^\alpha)\right] = g(t^\alpha)E_\alpha(-at^\alpha)$$

In the above steps we have used $^J_0 D^\alpha_t \left[E_\alpha(-at^\alpha)\right] = -aE_\alpha(-at^\alpha)$. Now operating $^J D^{-\alpha}$ on both the sides of the obtained last expression in above derivation i.e.

$^J D^\alpha \left[ y E_\alpha(-at^\alpha) \right] = g(t^\alpha) E_\alpha(-at^\alpha)$. Also we add a constant $A$ since Jumarie type derivative of a constant is zero and from here we get the following

$$\left[ y E_\alpha(-at^\alpha) \right] = {}^J D^{-\alpha} \left[ g(t^\alpha) E_\alpha(-at^\alpha) \right] + A \quad \text{where } A \text{ is a constant.}$$

$$y = \left[ E_\alpha(at^\alpha) \right] \left[ {}^J D^{-\alpha} \left( g(t^\alpha) E_\alpha(-at^\alpha) \right) + A \right] \tag{3.2}$$

or

$$y = A E_\alpha(at^\alpha) + \left[ E_\alpha(at^\alpha) \right] \left[ {}^J D^{-\alpha} \left( g(t^\alpha) E_\alpha(-at^\alpha) \right) \right]$$

the first part corresponds to solution of corresponding homogeneous equation, that is $y_c = A E_\alpha(at^\alpha)$ and the other part $y_p = \left[ E_\alpha(at^\alpha) \right] \left[ {}^J D^{-\alpha} \left( g(t^\alpha) E_\alpha(-at^\alpha) \right) \right]$ corresponds to the effect of non-homogeneous part and free from integral constant, but depending on the nature of forcing function, this part is named as Particular Integral (PI) as in case of classical differential equations. Now we take several forms of forcing function.

### 3.10 Particular Integral for $g(t^\alpha) = E_\alpha(ct^\alpha)$

Here consider the linear first order non-homogeneous fractional differential equation of order $\alpha$ with $0 < \alpha \leq 1$ with $y = 0$ for $t < 0$

$$({}^J D^\alpha - a) y = g(t^\alpha) \qquad g(t^\alpha) = E_\alpha(ct^\alpha), \qquad c \neq a.$$

then the Particular Integral (PI) described in the previous section is

$$y_p = \left[ E_\alpha(at^\alpha) \right] \left[ {}^J D^{-\alpha} \left( g(t^\alpha) E_\alpha(-at^\alpha) \right) \right]$$

Putting $g(t^\alpha) = E_\alpha(ct^\alpha)$ in above we get the following

$$\begin{aligned} y_p &= \left[ E_\alpha(at^\alpha) \right] \left[ {}^J D^{-\alpha} \left( E_\alpha(ct^\alpha) E_\alpha(-at^\alpha) \right) \right] \\ &= \left[ E_\alpha(at^\alpha) \right] \left[ {}^J D^{-\alpha} \left( E_\alpha((c-a)t^\alpha) \right) \right] \\ &= \left[ E_\alpha(at^\alpha) \right] \left( \frac{1}{c-a} \right) \left[ E_\alpha((c-a)t^\alpha) \right] \\ &= \frac{1}{c-a} E_\alpha(ct^\alpha) \end{aligned}$$

For $c = a$, P.I. is

$$y_p = \left[E_\alpha(at^\alpha)\right]\left[{}^JD^{-\alpha}\left(E_\alpha(at^\alpha)E_\alpha(-at^\alpha)\right)\right]$$

$$= \left[E_\alpha(at^\alpha)\right]\left({}^JD^{-\alpha}[1]\right) \quad \text{use} \quad {}^JD^{-\alpha}[1] = \frac{t^\alpha}{\Gamma(1+\alpha)}$$

$$= \frac{t^\alpha}{\Gamma(1+\alpha)} E_\alpha(at^\alpha)$$

- **Short procedure for calculating Particular Integral for** $g(t^\alpha) = E_\alpha(ct^\alpha)$.

This procedure is similar and in conjugation with classical integer order calculus. In classical order calculus $\alpha = 1$. Hence the forced function reduce to $g(t) = \exp(ct)$. Therefore the particular integral will be

$$y_p = \exp(at) \times \left[D^{-1}\left(\exp(ct).\exp(-at)\right)\right] = \left[\exp(at)\right]\left[\frac{1}{D}\left(\exp((c-a)t)\right)\right]$$

$$= \left[\exp(at)\right]\left(\frac{1}{c-a}\right)\left[\exp((c-a)t)\right] = \left(\frac{1}{c-a}\right)\left[\exp(ct)\right] \quad \text{for } c \neq a.$$

For $c = a$

$$y_p = \left[\exp(at)\right]\left[\frac{1}{D}\left(\exp(at)\exp(-at)\right)\right] = \left[\exp(at)\right]\left[\frac{1}{D}(1)\right] = t\exp(at)$$

Here we observe that the derivative operator $D$ is replaced by $c$ in the first case, i.e. for $c \neq a$. In the second case the derivative operator $D$ is replaced by $D + a$. We can replace the fractional Jumarie derivative operator ${}^JD^\alpha$ by $c$ for the first case $c \neq a$ and by ${}^JD^\alpha + a$ for second case $c = a$ The short procedure as follows for Particular Integral that is,

For $c \neq a$, $\quad y_p = \left(\dfrac{1}{{}^JD^\alpha - a}\right)\left[E_\alpha(ct^\alpha)\right] \quad$ replace ${}^JD^\alpha$ by $c$

$$= \frac{1}{c-a} E_\alpha(ct^\alpha)$$

For $c = a$, $\quad y_p = \left(\dfrac{1}{{}^JD^\alpha - a}\right)\left[E_\alpha(at^\alpha)\right] \quad$ replace ${}^JD^\alpha$ by ${}^JD^\alpha + a$

$$= \left[E_\alpha(at^\alpha)\right]\left(\frac{1}{{}^JD^\alpha + a - a}[1]\right)$$

$$= \left[E_\alpha(at^\alpha)\right]\left(\frac{1}{{}^JD^\alpha}[1]\right)$$

$$= \left[E_\alpha(at^\alpha)\right]\left({}^JD^{-\alpha}[1]\right)$$

$$= \frac{t^\alpha}{\Gamma(1+\alpha)} E_\alpha(at^\alpha).$$

Hence the general solution of equation (3.1) is $y = y_c + y_p$

$$y = \begin{cases} AE_\alpha(at^\alpha) + \dfrac{1}{c-a} E_\alpha(ct^\alpha) & \text{for } c \neq a \\ AE_\alpha(at^\alpha) + \dfrac{t^\alpha}{\Gamma(1+\alpha)} E_\alpha(at^\alpha) & \text{for } c = a. \end{cases}$$

**3.11 Particular Integral for** $g(t^\alpha) = t^\alpha$

Again when $g(t^\alpha) = t^\alpha$ then the differential equation (3.1) becomes

$$({}^J D^\alpha - a) y = t^\alpha \qquad (3.3)$$

The solution of the homogeneous part [15] that is $({}^J D^\alpha - a) y = 0$ is $y_c = A_1 E_\alpha(at^\alpha)$

Let $y = V(t^\alpha) E_\alpha(at^\alpha)$ the solution of the corresponding non-homogeneous equation where $V(t^\alpha)$ is an unknown function of $t^\alpha$. Then using the definition by Jumarie [16] that is

$$(u(x)v(x))^{(\alpha)} = u^{(\alpha)}(x) v(x) + u(x) v^{(\alpha)}(x)$$

We get the following

$$\begin{aligned}{}^J D^\alpha[y] &= {}^J D^\alpha \left[ V(t^\alpha) E_\alpha(at^\alpha) \right] \\ &= \left[ V(t^\alpha) \right] \left( {}^J D^\alpha \left[ E_\alpha(at^\alpha) \right] \right) + \left[ E_\alpha(at^\alpha) \right] \left( {}^J D^\alpha \left[ V(t^\alpha) \right] \right) \\ &= \left[ V(t^\alpha) \right] \left[ aE_\alpha(at^\alpha) \right] + \left[ E_\alpha(at^\alpha) \right] \left( {}^J D^\alpha \left[ V(t^\alpha) \right] \right) \end{aligned}$$

putting this in (3.3) we get

$$({}^J D^\alpha - a) y = t^\alpha$$
$${}^J D^\alpha[y] - ay = t^\alpha$$
$$\left[ V(t^\alpha) \right] \left[ aE_\alpha(at^\alpha) \right] + \left[ E_\alpha(at^\alpha) \right] \left( {}^J D^\alpha \left[ V(t^\alpha) \right] \right) - a \left[ V(t^\alpha) \right] \left[ E_\alpha(at^\alpha) \right] = t^\alpha$$
$$\left[ E_\alpha(at^\alpha) \right] \left( {}^J D^\alpha \left[ V(t^\alpha) \right] \right) = t^\alpha$$

Therefore we get

$$ {}^J D^\alpha \left[ V(t^\alpha) \right] = \frac{t^\alpha}{E_\alpha(at^\alpha)} = t^\alpha E_\alpha(-at^\alpha)$$

We now apply fractional integration by parts by Jumarie formula [16] as depicted below

$$\int_0^x u(y)\,{}^J D^\alpha[v(y)](dy)^\alpha = u(y)v(y)\big|_0^x - \int_0^x v(y)\,{}^J D^\alpha[u(y)](dy)^\alpha$$

Here we mention that the symbol $\int f(x)(dx)^\alpha \equiv {}^J D^{-\alpha} f(x)$ implies Jumarie fractional integration as defined in section-2. We will use also [15] derived expression that is ${}^J D^{-\alpha}\left[E_\alpha(-at^\alpha)\right] = (-a^{-1})\left[E_\alpha(-at^\alpha)\right]$, in the following derivation.

$$V(t^\alpha) = A + {}^J D^{-\alpha}\left[t^\alpha E_\alpha(-at^\alpha)\right]$$

$$= A + \int_0^t \tau^\alpha E_\alpha(-a\tau^\alpha)(d\tau)^\alpha$$

$$= A + \int_0^t (\tau^\alpha)\left({}^J D^\alpha\left[\tfrac{E_\alpha(-a\tau^\alpha)}{-a}\right]\right)(d\tau)^\alpha$$

$$= A + \tau^\alpha \frac{E_\alpha(-a\tau^\alpha)}{-a}\Big|_0^t - \int_0^t \frac{E_\alpha(-a\tau^\alpha)}{-a}\left({}^J D_\tau^\alpha\left[\tau^\alpha\right]\right)(d\tau)^\alpha \quad \text{use} \quad {}^J D_x^\alpha\left[x^\alpha\right] = \Gamma(1+\alpha)$$

$$= A + \tau^\alpha \frac{E_\alpha(-a\tau^\alpha)}{-a}\Big|_0^t - \int_0^t \frac{E_\alpha(-a\tau^\alpha)}{-a}\Gamma(1+\alpha)(d\tau)^\alpha$$

$$= A + \tau^\alpha \frac{E_\alpha(-a\tau^\alpha)}{-a}\Big|_0^t + \frac{\Gamma(1+\alpha)}{a}\left[\frac{E_\alpha(-a\tau^\alpha)}{-a}\right]_0^t$$

$$= A - t^\alpha \frac{E_\alpha(-at^\alpha)}{a} - \frac{\Gamma(1+\alpha)}{a^2}\left[E_\alpha(-at^\alpha) - 1\right]$$

$A =$ constant.

Hence the general solution is

$$y = \left[V(t^\alpha)\right]\left[E_\alpha(at^\alpha)\right]$$

$$= -\frac{1}{a}\left(t^\alpha + \frac{\Gamma(1+\alpha)}{a}\right) + \left(A + \frac{\Gamma(1+\alpha)}{a^2}\right)E_\alpha(at^\alpha)$$

$$= A_1 E_\alpha(at^\alpha) - \frac{1}{a}\left(t^\alpha + \frac{\Gamma(1+\alpha)}{a}\right) \quad \text{where} \quad A_1 = A + \frac{\Gamma(1+\alpha)}{a^2}$$

$y_c = A_1 E_\alpha(at^\alpha)$ the first part in above expression is solution of homogeneous equation and the second part of the above that is $y_p = -\frac{1}{a}\left(t^\alpha + \frac{\Gamma(1+\alpha)}{a}\right)$ is particular integral.

- **Short procedure for Calculating Particular Integral for** $g(t^\alpha) = t^\alpha$

This procedure is similar and in conjugation with classical integer order calculus. Here for $\alpha = 1$, and $g(t) = t$, and the corresponding particular integral is

$$y_p = \frac{1}{D-a}t = -\frac{1}{a}\left(1-\frac{D}{a}\right)^{-1}t = -\frac{1}{a}\left\{1+\frac{D}{a}+\frac{D^2}{a^2}+...\right\}t = -\frac{1}{a}\left(t+\frac{1}{a}\right)$$

In the same way we can have a short procedure as follows for Particular Integral that is,

$$y_p = \frac{1}{{}^J D^\alpha - a}t^\alpha$$

$$= -\frac{1}{a}\left(1-\frac{{}^J D^\alpha}{a}\right)^{-1}t^\alpha$$

$$= -\frac{1}{a}\left(1+\frac{{}^J D^\alpha}{a}+\frac{{}^J D^{2\alpha}}{a^2}+...\right)t^\alpha$$

$$= -\frac{1}{a}\left(t^\alpha + \frac{\Gamma(1+\alpha)}{a}\right)$$

In the above derivation ${}^J D^{2\alpha}[t^\alpha] = {}^J D^\alpha[D^\alpha(t^\alpha)] = {}^J D^\alpha[\Gamma(1+\alpha)] = D^\alpha[C] = 0$ is used. Thus all the Jumarie derivatives ${}^J D^{k\alpha}[t^\alpha] = 0$ for $k > 1$, where $k$ is Natural number. Therefore we have discussed the solutions of non-homogeneous $\alpha$-order differential equations for different forcing functions $g(t^\alpha)$.

**3.12 Evaluation of** $y_p = \frac{1}{{}^J D^{2\alpha} - a^2}\sin_\alpha(ct^\alpha)$, $g(t^\alpha) = \sin_\alpha(ct^\alpha)$ where $c^2 + a^2 \neq 0$

$\frac{1}{{}^J D^{2\alpha} - a^2}$ can be factorized as $\frac{1}{{}^J D^{2\alpha} - a^2} = \frac{1}{({}^J D^\alpha - a)({}^J D^\alpha + a)} = \frac{1}{2a}\left[\frac{1}{{}^J D^\alpha - a} - \frac{1}{{}^J D^\alpha + a}\right]$, and we use this in following derivation.

As in section 3.10 here we replace ${}^J D^\alpha$ by $ic$ and by $-ic$ for the operations $\frac{1}{{}^J D^\alpha - a}E_\alpha(ict^\alpha)$ and $\frac{1}{{}^J D^\alpha - a}E_\alpha(-ict^\alpha)$ respectively, as is demonstrated below.

$$y_p = \frac{1}{{}^J D^{2\alpha} - a^2}\sin_\alpha(ct^\alpha)$$

$$= \frac{1}{2a}\left[\frac{1}{{}^J D^\alpha - a} - \frac{1}{{}^J D^\alpha + a}\right]\sin_\alpha(ct^\alpha) = \frac{1}{2a}\left[\frac{1}{{}^J D^\alpha - a}\sin_\alpha(ct^\alpha) - \frac{1}{{}^J D^\alpha + a}\sin_\alpha(ct^\alpha)\right]$$

$$\sin_\alpha(ct^\alpha) \stackrel{def}{=} \frac{1}{2i}\left[E_\alpha(ict^\alpha) - E_\alpha(-ict^\alpha)\right]$$

$$\frac{1}{{}^J D^\alpha - a} \sin_\alpha(ct^\alpha) = \frac{1}{2i}\left[\frac{1}{{}^J D^\alpha - a} E_\alpha(ict^\alpha) - \frac{1}{{}^J D^\alpha - a} E_\alpha(-ict^\alpha)\right]$$

$$= \frac{1}{2i}\left[\frac{1}{ic - a} E_\alpha(ict^\alpha) - \frac{1}{-ic - a} E_\alpha(-ict^\alpha)\right]$$

$$\frac{1}{{}^J D^\alpha + a} \sin_\alpha(ct^\alpha) = \frac{1}{2i}\left[\frac{1}{{}^J D^\alpha + a} E_\alpha(ict^\alpha) - \frac{1}{{}^J D^\alpha + a} E_\alpha(-ict^\alpha)\right]$$

$$= \frac{1}{2i}\left[\frac{1}{ic + a} E_\alpha(ict^\alpha) - \frac{1}{-ic + a} E_\alpha(-ict^\alpha)\right]$$

Therefore

$$\frac{1}{{}^J D^{2\alpha} - a^2} \sin_\alpha(ct^\alpha)$$

$$= \frac{1}{4ai}\left[\frac{1}{ic - a} E_\alpha(it^\alpha) - \frac{1}{ic + a} E_\alpha(ict^\alpha) - \frac{1}{-ic - a} E_\alpha(-ict^\alpha) + \frac{1}{-ic + a} E_\alpha(-ict^\alpha)\right]$$

$$= \frac{1}{4ai}\left[\left(\frac{1}{ic - a} - \frac{1}{ic + a}\right) E_\alpha(ict^\alpha) - \left(\frac{1}{-ic - a} - \frac{1}{-ic + a}\right) E_\alpha(-ict^\alpha)\right]$$

$$= \frac{1}{2i}\left[\frac{1}{-c^2 - a^2} E_\alpha(ict^\alpha) - \frac{1}{-c^2 - a^2} E_\alpha(-ict^\alpha)\right] = \left(\frac{1}{-c^2 - a^2}\right)\left(\frac{1}{2i}\left[E_\alpha(ict^\alpha) - E_\alpha(-ict^\alpha)\right]\right)$$

$$= \frac{\sin_\alpha(ct^\alpha)}{-c^2 - a^2}$$

Similarly we get by following above procedure

$$\frac{1}{{}^J D^{2\alpha} - a^2} \cos_\alpha(ct^\alpha) = \frac{\cos_\alpha(ct^\alpha)}{-c^2 - a^2}$$

Thus to find the particular integral $\frac{1}{{}^J D^{2\alpha} - a^2} \sin_\alpha(ct^\alpha)$ replace ${}^J D^{2\alpha}$ by $-c^2$.

This procedure is similar and in conjugation with classical integer order calculus. In classical order calculus $\alpha = 1$ hence the forced function reduce to $g(t) = \sin(ct)$. Therefore the particular integral will be

$$y_p = \frac{1}{{}^J D^2 - a^2} \sin(ct)$$

$$= \frac{1}{2a} \left[ \frac{1}{{}^J D - a} - \frac{1}{{}^J D + a} \right] \sin(ct) = \frac{1}{2a} \left[ \frac{1}{{}^J D - a} \sin(ct) - \frac{1}{{}^J D + a} \sin(ct) \right]$$

$$\sin(ct) = \tfrac{1}{2i} \left[ \exp(ict) - \exp(-ict) \right]$$

$$\frac{1}{{}^J D - a} \sin(ct) = \frac{1}{2i} \left[ \frac{1}{{}^J D - a} \exp(ict) - \frac{1}{{}^J D - a} \exp(-ict) \right]$$

$$= \frac{1}{2i} \left[ \frac{1}{ic - a} \exp(ict) - \frac{1}{-ic - a} \exp(-ict) \right]$$

$$\frac{1}{{}^J D + a} \sin(ct) = \frac{1}{2i} \left[ \frac{1}{{}^J D + a} \exp(ict) - \frac{1}{{}^J D + a} \exp(-ict) \right]$$

$$= \frac{1}{2i} \left[ \frac{1}{ic + a} \exp(ict) - \frac{1}{-ic + a} \exp(-ict) \right]$$

Therefore

$$\frac{1}{{}^J D^2 - a^2} \sin(ct)$$

$$= \frac{1}{4ai} \left[ \frac{1}{ic - a} \exp(ict) - \frac{1}{ic + a} \exp(ict) - \frac{1}{-ic - a} \exp(-ict) + \frac{1}{-ic + a} \exp(-ict) \right]$$

$$= \frac{1}{4ai} \left[ \left( \frac{1}{ic - a} - \frac{1}{ic + a} \right) \exp(ict) - \left( \frac{1}{-ic - a} - \frac{1}{-ic + a} \right) \exp(ict) \right]$$

$$= \frac{1}{2i} \left[ \frac{1}{-c^2 - a^2} \exp(ict) - \frac{1}{-c^2 - a^2} \exp(-ict) \right] = \left( \frac{1}{-c^2 - a^2} \right) \left( \frac{1}{2i} [\exp(ict) - \exp(-ict)] \right)$$

$$= \frac{\sin(ct)}{-c^2 - a^2}$$

## 4.0 $2\alpha-$ order non-homogeneous fractional differential equations

General formulation of non-homogeneous fractional differential equation of

$2\alpha -$ order is $\quad {}^J D^{2\alpha}[y] + p({}^J D^{\alpha}[y]) + qy = g(t^{\alpha}) \quad 0 < \alpha < 1 \quad y(t) = 0 \quad \text{for} \quad t < 0$

where $p$ and $q$ are constant here. Consider the $2\alpha-$ order non-homogeneous fractional differential equation $f({}^J D^{\alpha}) y = g(t^{\alpha})$ where $f({}^J D^{\alpha}) = ({}^J D^{\alpha} - a)({}^J D^{\alpha} - b)$ then solution of the non-homogeneous part that is $f({}^J D^{\alpha}) = 0$ given by $y_c = AE_{\alpha}(at^{\alpha}) + BE_{\alpha}(bt^{\alpha})$ [15].

## 4.10 Use of method of un-determinant coefficient method to calculate the particular integrals for different functional forms of $g(t^\alpha)$.

For $g(t^\alpha) = E_\alpha(ct^\alpha)$ we have the given equation is

$$f(^JD^\alpha)y = (^JD^\alpha - a)(^JD^\alpha - b)y$$
$$= \left[^JD^{2\alpha} - (a+b)(^JD^\alpha) + ab\right](y) \quad (4.1)$$
$$= E_\alpha(ct^\alpha) \quad \text{for } c \neq a \neq b$$

Here let the particular integral be $y_p = PE_\alpha(ct^\alpha)$ where $P$ is constant.

Then

$$^JD^\alpha\left[y_p\right] = PcE_\alpha(ct^\alpha) \qquad ^JD^{2\alpha}\left[y_p\right] = Pc^2 E_\alpha(ct^\alpha)$$

and putting in the given equation (4.1) we get the following

$$\left[^JD^{2\alpha} - (a+b)(^JD^\alpha) + ab\right](y_p) = E_\alpha(ct^\alpha)$$
$$^JD^{2\alpha}\left[y_p\right] - (a+b)(^JD^\alpha\left[y_p\right]) + aby_p = E_\alpha(ct^\alpha)$$
$$Pc^2 E_\alpha(ct^\alpha) - (a+b)PcE_\alpha(ct^\alpha) + abPE_\alpha(ct^\alpha) = E_\alpha(ct^\alpha)$$
$$Pc^2 - (a+b)Pc + abP = 1$$
$$P\left[c^2 - (a+b)c + ab\right] = 1$$
$$P\left[(c-a)(c-b)\right] = 1$$

Therefore

$$P = \frac{1}{(c-a)(c-b)}$$

and consequently the Particular integral is

$$y_p = \frac{1}{(c-a)(c-b)} E_\alpha(ct^\alpha)$$

Hence the general solution is

$$y = AE_\alpha(at^\alpha) + BE_\alpha(bt^\alpha) + \frac{1}{(c-a)(c-b)} E_\alpha(ct^\alpha) \quad (4.2)$$

For $c = a \neq b$ implying $(c-a)(c-b) = 0$ then the solution (4.2) does not exists. In this case the fractional differential equation is

$$f(^JD^\alpha)y = (^JD^\alpha - a)(^JD^\alpha - b)y = \left[^JD^{2\alpha} - (a+b)(^JD^\alpha) + ab\right]y = E_\alpha(at^\alpha) \tag{4.3}$$

If we consider $y_p = PE_\alpha(at^\alpha)$ in this case also then putting in (4.3) we get $0 = E_\alpha(at^\alpha)$ which is free from $P$ i.e. $P$ is non-determinable. This form of PI is not suitable here; consider the modified form as following

$$y_p = Pt^\alpha E_\alpha(at^\alpha)$$

Then

$$^JD^\alpha[y_p] = P\left[at^\alpha E_\alpha(at^\alpha) + \Gamma(1+\alpha)E_\alpha(at^\alpha)\right]$$
$$^JD^{2\alpha}[y_p] = P\left[a^2 t^\alpha E_\alpha(at^\alpha) + 2a\Gamma(1+\alpha)E_\alpha(at^\alpha)\right]$$

and putting in (4.3) we get the following

$$Pa^2 t^\alpha + 2aP\Gamma(1+\alpha) - (a+b)Pat^\alpha - (a+b)P\Gamma(1+\alpha) + Pabt^\alpha = 1$$
$$P = \frac{1}{(a-b)\Gamma(1+\alpha)}$$
$$y_p = \frac{t^\alpha}{(a-b)\Gamma(1+\alpha)} E_\alpha(at^\alpha)$$

In this case the general solution of the fractional differential equation is of following form

$$y = AE_\alpha(at^\alpha) + BE_\alpha(bt^\alpha) + \frac{t^\alpha}{(a-b)\Gamma(1+\alpha)} E_\alpha(at^\alpha)$$

When $c = a = b$ then

$$f(^JD^\alpha)y = (^JD^\alpha - a)^2 y = E_\alpha(at^\alpha) \tag{4.4}$$

and take the particular integral in the form

$$y_p = Pt^{2\alpha} E_\alpha(at^\alpha)$$

then

$$D^\alpha y_p = Pat^{2\alpha} E_\alpha(at^\alpha) + Pt^\alpha \frac{\Gamma(1+2\alpha)}{\Gamma(1+\alpha)} E_\alpha(at^\alpha)$$

$$^JD^{2\alpha}[y_p] = Pa^2 t^{2\alpha} E_\alpha(at^\alpha) + 2aPt^\alpha \frac{\Gamma(1+2\alpha)}{\Gamma(1+\alpha)} E_\alpha(at^\alpha) + P\Gamma(1+2\alpha)E_\alpha(at^\alpha)$$

Putting this in (4.4) and after simplification we get

$$P = \frac{1}{\Gamma(1+2\alpha)}$$

$$y_p = \frac{t^{2\alpha}}{\Gamma(1+2\alpha)} E_\alpha(at^\alpha)$$

In this case the general solution will be of following form

$$y = (At^\alpha + B)E_\alpha(at^\alpha) + \frac{t^{2\alpha}}{\Gamma(1+2\alpha)} E_\alpha(at^\alpha)$$

Thus we can summarize the result as a theorem in the following form

**Theorem:** The differential equation the $f(^JD^\alpha)y = E_\alpha(ct^\alpha)$ has particular integral

$$\frac{1}{f(c)} E_\alpha(ct^\alpha) \quad \text{for} \quad f(c) \neq 0$$

(i) When $f(^JD^\alpha) = (^JD^\alpha - a)(^JD^\alpha - b)$ for $c \neq a \neq b$ then solution of the fractional differential equation will be

$$y = AE_\alpha(at^\alpha) + BE_\alpha(bt^\alpha) + \frac{1}{(c-a)(c-b)} E_\alpha(ct^\alpha)$$

(ii) When $f(^JD^\alpha) = (^JD^\alpha - a)(^JD^\alpha - b)$ for $c = a \neq b$ then solution of the fractional differential equation will be

$$y = AE_\alpha(at^\alpha) + BE_\alpha(bt^\alpha) + \frac{t^\alpha}{(a-b)\Gamma(1+\alpha)} E_\alpha(at^\alpha)$$

(iii) When $f(^JD^\alpha) = (^JD^\alpha - a)(^JD^\alpha - b)$ for $c = a = b$ then solution of the fractional differential equation will be

$$y = (At^\alpha + B)E_\alpha(at^\alpha) + \frac{t^{2\alpha}}{\Gamma(1+2\alpha)} E_\alpha(at^\alpha)$$

**4.11 Use of direct method to calculate the Particular integrals for different functional format of $g(t^\alpha)$.**

Using the direct method as describe in section 3.10 we can easily calculate the Particular integrals for different functional format of $g(t^\alpha)$.

- For $g(t^\alpha) = E_\alpha(ct^\alpha)$ we have

$$y_p = \frac{1}{f({}^J D^\alpha)} E_\alpha(ct^\alpha)$$

$$= \frac{1}{({}^J D^\alpha - a)({}^J D^\alpha - b)} E_\alpha(ct^\alpha)$$

$$= \frac{1}{a-b}\left[\frac{1}{({}^J D^\alpha - a)} - \frac{1}{({}^J D^\alpha - b)}\right] E_\alpha(ct^\alpha)$$

$$= \frac{1}{a-b}\left[\frac{1}{(c-a)} - \frac{1}{(c-b)}\right] E_\alpha(ct^\alpha) \qquad \text{for } c \neq a \neq b$$

$$= \frac{1}{(c-b)(c-a)} E_\alpha(ct^\alpha) \qquad \text{for } c \neq a \neq b.$$

In this case the general solution is

$$y = AE_\alpha(at^\alpha) + BE_\alpha(bt^\alpha) + \frac{1}{(c-b)(c-a)} E_\alpha(ct^\alpha) \qquad \text{for } c \neq a$$

- For $c = a$

$$y_p = \frac{1}{f({}^J D^\alpha)} E_\alpha(at^\alpha)$$

$$= \frac{1}{({}^J D^\alpha - a)({}^J D^\alpha - b)} E_\alpha(at^\alpha)$$

$$= \frac{1}{a-b}\left[\frac{1}{({}^J D^\alpha - a)} - \frac{1}{({}^J D^\alpha - b)}\right] E_\alpha(at^\alpha)$$

$$= \frac{1}{a-b}\left[\frac{t^\alpha E_\alpha(at^\alpha)}{\Gamma(1+\alpha)} - \frac{1}{(a-b)}\right] E_\alpha(at^\alpha)$$

In this case the second part will be adjusted in the complementary function and hence the general solution is

$$y = AE_\alpha(at^\alpha) + BE_\alpha(bt^\alpha) + \frac{t^\alpha}{(a-b)\Gamma(1+\alpha)} E_\alpha(at^\alpha)$$

- For $c = a = b$

$$y_p = \frac{1}{f({}^JD^\alpha)} E_\alpha(at^\alpha)$$

$$= \frac{1}{\left({}^JD^\alpha - a\right)^2} E_\alpha(at^\alpha)$$

$$= \left[\frac{1}{\left({}^JD^\alpha - a\right)^2}\right] E_\alpha(at^\alpha)$$

$$= E_\alpha(at^\alpha) \frac{1}{\left({}^JD^\alpha + a - a\right)^2} [1]$$

$$= \frac{t^{2\alpha} E_\alpha(at^\alpha)}{\Gamma(1+2\alpha)}.$$

In this case the general solution will be

$$y = (At^\alpha + B)E_\alpha(at^\alpha) + \frac{t^{2\alpha}}{\Gamma(1+2\alpha)} E_\alpha(at^\alpha)$$

- In generalized case for any polynomial type function $f({}^JD^\alpha)$ the particular integral is

$$y_p = \frac{1}{f(D^\alpha)} E_\alpha(ct^\alpha) = \frac{1}{f(c)} E_\alpha(ct^\alpha) \quad \text{for} \quad f(c) \neq 0.$$

For $f(c) = 0$ $\quad f({}^JD^\alpha)$ must contain a factor of the form $\left({}^JD^\alpha - c\right)^k$, $k$-positive integer i.e.
$f({}^JD^\alpha) = \left({}^JD^\alpha - c\right)^k \phi({}^JD^\alpha)$ with $\phi(c) \neq 0$ then

$$y_p = \frac{1}{({}^JD^\alpha - c)^k \phi({}^JD^\alpha)} E_\alpha(ct^\alpha)$$

$$= \frac{E_\alpha(ct^\alpha)}{\phi(c)} \frac{1}{\left({}^JD^\alpha - c\right)^k} [1]$$

$$= \frac{E_\alpha(ct^\alpha)}{\phi(c)} \frac{t^{k\alpha}}{\Gamma(1+k\alpha)}$$

- Again if
  $f({}^JD^\alpha) = v(t^\alpha) E_\alpha(ct^\alpha)$ Here using Leibnitz rule of fractional derivative on (Jumarie type fractional derivative) we get the following steps

$$^J D^\alpha \left( v(t^\alpha) E_\alpha(ct^\alpha) \right) = \left( ^J D^\alpha \left[ v(t^\alpha) \right] \right) \left[ E_\alpha(ct^\alpha) \right] + cv(t^\alpha) E_\alpha(ct^\alpha)$$

$$= (^J D^\alpha + c) v(t^\alpha) E_\alpha(ct^\alpha)$$

operate $^J D^{-\alpha}$ on both sides

$$v(t^\alpha) E_\alpha(ct^\alpha) = {}^J D^{-\alpha} \left[ (^J D^\alpha + c) v(t^\alpha) E_\alpha(ct^\alpha) \right]$$

Let $\quad V(t^\alpha) = (^J D^\alpha + c) v(t^\alpha)$

then $\quad v(t^\alpha) = \dfrac{1}{^J D^\alpha + c} V(t^\alpha)$

$$^J D^{-\alpha} \left[ V(t^\alpha) E_\alpha(ct^\alpha) \right] = \dfrac{E_\alpha(ct^\alpha)}{^J D^\alpha + c} V(t^\alpha)$$

Thus we obtain

$$^J D^{-\alpha} \left[ V(t^\alpha) E_\alpha(ct^\alpha) \right] = \dfrac{E_\alpha(ct^\alpha)}{D^\alpha + c} V(t^\alpha)$$

$$\dfrac{1}{^J D^\alpha} V(t^\alpha) E_\alpha(ct^\alpha) = \dfrac{E_\alpha(ct^\alpha)}{^J D^\alpha + c} V(t^\alpha)$$

Thus in generalized case we have

$$\dfrac{1}{f(^J D^\alpha)} V(t^\alpha) E_\alpha(ct^\alpha) = \dfrac{E_\alpha(ct^\alpha)}{f(^J D^\alpha + c)} V(t^\alpha)$$

## 5.0 Solution the fractional differential-application of method derived

**Example 1**: We take the following fractional differential equation $y(t) = 0$ for $t < 0$

$$^J D^{2/3} [y] - 5 \left( D^{1/3} [y] \right) + 6 y = t^2$$

**Solution**: Solution of the corresponding homogeneous equation is [15]

$$y_c = AE_{1/3}(2t^{1/3}) + BE_{1/3}(3t^{1/3})$$

The particular integral calculation is done in following steps

$$y_p = \frac{1}{(^J D^{1/3} - 2)(^J D^{1/3} - 3)} t^{6/3}$$

$$= \left[ \frac{1}{(^J D^{1/3} - 3)} - \frac{1}{(^J D^{1/3} - 2)} \right] t^{6/3}$$

$$= \left[ -\frac{1}{3}\left(1 - \frac{^J D^{1/3}}{3}\right)^{-1} + \frac{1}{2}\left(1 - \frac{^J D^{1/3}}{2}\right)^{-1} \right] t^{6/3}$$

$$= -\frac{1}{3}\left(1 - \frac{^J D^{1/3}}{3} + \frac{^J D^{2/3}}{9} - \frac{^J D^{3/3}}{27} + \ldots\right) t^{6/3} + \frac{1}{2}\left(1 - \frac{^J D^{1/3}}{2} + \frac{^J D^{2/3}}{4} - \frac{^J D^{3/3}}{8} + \ldots\right) t^{6/3}$$

$$= \frac{1}{6} t^{6/3} + \frac{3^2 - 2^2}{3^2 \cdot 2^2} \frac{\Gamma(1+2)}{\Gamma\left(1+\frac{5}{3}\right)} t^{5/3} - \frac{3^3 - 2^3}{3^3 \cdot 2^3} \frac{\Gamma(1+2)}{\Gamma\left(1+\frac{4}{3}\right)} t^{4/3} + \frac{3^4 - 2^4}{3^4 \cdot 2^4} \frac{\Gamma(1+2)}{\Gamma\left(1+\frac{3}{3}\right)} t^{3/3} +$$

$$+ \frac{3^5 - 2^5}{3^5 \cdot 2^5} \frac{\Gamma(1+2)}{\Gamma\left(1+\frac{2}{3}\right)} t^{2/3} + \frac{3^6 - 2^6}{3^6 \cdot 2^6} \frac{\Gamma(1+2)}{\Gamma\left(1+\frac{1}{3}\right)} t^{1/3} + \frac{3^7 - 2^7}{3^7 \cdot 2^7} \frac{\Gamma(1+2)}{\Gamma(1)}$$

Hence the general solution is

$$y = A E_{1/3}(2 t^{1/3}) + B E_{1/3}(3 t^{1/3}) + \frac{1}{6} t^{6/3} + \frac{10}{6^2 \Gamma\left(\frac{8}{3}\right)} t^{5/3} - \frac{38}{6^3 \Gamma\left(\frac{7}{3}\right)} t^{4/3} + \frac{130}{6^4 \Gamma\left(\frac{6}{3}\right)} t^{3/3} +$$

$$+ \frac{422}{6^5 \Gamma\left(\frac{5}{3}\right)} t^{2/3} + \frac{1330}{6^6 \Gamma\left(\frac{4}{3}\right)} t^{1/3} + \frac{4118}{6^7}$$

Where $A$ and $B$ are arbitary constants.

**Example 2**: Consider the fractional order forced differential equation $y(t) = 0$ for $t < 0$

$$^J D^{2\alpha}[y] + \omega^2 y = F \cos_\alpha(a t^\alpha) \quad a \neq \omega$$

**Solution**: Here solution of the corresponding homogeneous equation $^J D^{2\alpha}[y] + \omega^2 y = 0$ is [15]

$$y_c = A \cos_\alpha(\omega t^\alpha) + B \sin_\alpha(\omega t^\alpha)$$

The particular integral is

$$y_p = \frac{1}{D^{2\alpha} + \omega^2} F \cos_\alpha(a t^\alpha)$$

$$= \frac{1}{-a^2 + \omega^2} F \cos_\alpha(a t^\alpha)$$

Hence the general solution is

$$y = A\cos_\alpha(\omega t^\alpha) + B\sin_\alpha(\omega t^\alpha) + \frac{1}{-a^2 + \omega^2} F\cos_\alpha(at^\alpha)$$

For particular integral we replaced $^JD^\alpha \equiv ia$ so $^JD^{2\alpha} \equiv (ia)^2 = -a^2$

**Example 3**: Take the fractional order damped-forced differential equation $y(t) = 0$ for $t < 0$

$$^JD^{2\alpha}[y] + 2c(D^\alpha[y]) + (c^2 + \omega^2)y = F\cos_\alpha(at^\alpha) \qquad a \neq \omega$$

**Solution**: Here solution of the corresponding homogeneous equation [15]

$$^JD^{2\alpha}[y] + 2c(D^\alpha[y]) + (c^2 + \omega^2)y = 0$$

is

$$y_c = [E_\alpha(ct^\alpha)][A\cos_\alpha(\omega t^\alpha) + B\sin_\alpha(\omega t^\alpha)]$$

The particular integral is

$$y_p = \frac{1}{^JD^{2\alpha} + 2c(^JD^\alpha) + (c^2 + \omega^2)} F\cos_\alpha(at^\alpha)$$

$$= \frac{1}{2c(^JD^\alpha) + (c^2 - a^2 + \omega^2)} F\cos_\alpha(at^\alpha) \quad \text{, here replace } ^JD^{2\alpha} \text{ by } -a^2$$

multiply the numerator and denominator by $(c^2 - a^2 + \omega^2) - 2c(^JD^\alpha)$

$$y_p = \frac{(c^2 - a^2 + \omega^2) - 2c(^JD^\alpha)}{(c^2 - a^2 + \omega^2)^2 - 4c^2(^JD^{2\alpha})} F\cos_\alpha(at^\alpha), \quad \text{here again replace } ^JD^{2\alpha} \text{ by } -a^2$$

$$= F\frac{(c^2 - a^2 + \omega^2)\cos_\alpha(at^\alpha) - 2ca\sin_\alpha(at^\alpha)}{(c^2 - a^2 + \omega^2)^2 + 4c^2a^2}$$

Hence the general solution is

$$y = (A\cos_\alpha(\omega t^\alpha) + B\sin_\alpha(\omega t^\alpha))[E_\alpha(ct^\alpha)] + F\frac{(c^2 - a^2 + \omega^2)\cos_\alpha(at^\alpha) - 2ca\sin_\alpha(at^\alpha)}{(c^2 - a^2 + \omega^2)^2 + 4c^2a^2}$$

## Conclusions

In this paper we have developed a method to solve the linear fractional non-homogeneous fractional differential equations, composed by Jumarie type derivative. The solutions are obtained here in terms of Mittag-Leffler function and fractional sine fractional cosine functions. Here we have proved via usage of Jumarie fractional derivative operator that for obtaining the particular integrals for several forcing functions scaled in function of variable $t^\alpha$ eases the method, and we obtain conjugation with classical method to solve classical non-homogeneous differential equations. The short cut rules, that are developed here in this paper to replace the operator $D^\alpha$ or operator $D^{2\alpha}$ as were used in classical calculus, gives ease in evaluating particular integrals. These techniques obtained herein this paper is remarkable to study fractional dynamic systems, and eases to get solution in terms of Mittag-Leffler, and fractional-trigonometric functions as in conjugation with exponential and normal trigonometric function for normal integer order calculus. Therefore this developed method is useful as it is having conjugation with the classical methods of solving non-homogeneous fractional linear differential equations composed via Jumarie fractional derivative, and is also useful in understanding physical systems described by FDE.

## Acknowledgement

Acknowledgments are to **Board of Research in Nuclear Science** (BRNS), Department of Atomic Energy Government of India for financial assistance received through BRNS research project no. 37(3)/14/46/2014-BRNS with BSC BRNS, title "Characterization of unreachable (Holderian) functions via Local Fractional Derivative and Deviation Function.